\documentclass[12pt,a4paper,reqno]{amsart}

\pdfoutput=1

\usepackage[utf8]{inputenc}
\usepackage[english]{babel}
\usepackage{amssymb}
\usepackage{amsmath}
\usepackage[colorlinks]{hyperref}
\usepackage{enumitem}
\usepackage{faktor}
\usepackage{mathtools}
\usepackage{stackrel}
\usepackage{verbatim}
\usepackage{tikz-cd}
\usepackage[english]{cleveref}
\usepackage{todonotes}
\usepackage{scalerel}
\usepackage{mathscinet}

\theoremstyle{definition}
\newtheorem{theorem}{Theorem}[section]

\newtheorem{definition}[theorem]{Definition}
\newtheorem{remark}[theorem]{Remark}

\newtheorem{lemma}[theorem]{Lemma}
\newtheorem{proposition}[theorem]{Proposition}
\newtheorem{corollary}[theorem]{Corollary}
\newtheorem{conjecture}[theorem]{Conjecture}

\newcommand{\Q}{\mathbb{Q}}
\newcommand{\C}{\mathbb{C}}
\newcommand{\R}{\mathbb{R}}
\newcommand{\dd}{\mathrm{d}}

\newcommand{\sgn}{\mathrm{sgn}}

\DeclareMathOperator{\Exp}{Exp}
\DeclareMathOperator{\ch}{ch}
\DeclareMathOperator{\Conf}{Conf} 
\DeclareMathOperator{\Ind}{Ind}
\DeclareMathOperator{\Res}{Res}

\allowdisplaybreaks

\makeatletter
\g@addto@macro{\UrlBreaks}{\UrlOrds}
\g@addto@macro{\UrlBreaks}{%
\do\/\do\d%
}
\makeatother

\begin{document}
\title[Orlik-Terao algebra of type A]{The Frobenius character of the Orlik-Terao algebra of type A}
\author[R. Pagaria]{Roberto Pagaria}
\thanks{The author is partially supported by PRIN 2017YRA3LK, by H2020 MCSA RISE project GHAIA -  n. 777822, and by National Science Foundation under Grant No. DMS-1929284 }
\address{Roberto Pagaria \newline \textup{Università di Bologna, Dipartimento di Matematica}\\ Piazza di Porta San Donato 5 - 40126 Bologna\\ Italy.}
\email{roberto.pagaria@unibo.it}

\begin{abstract}
We provide a new virtual description of the symmetric group action on the cohomology of ordered configuration space on $SU_2$ up to translations.
We use this formula to prove the Moseley-Proudfoot-Young conjecture.
As a consequence we obtain the graded Frobenius character of the Orlik-Terao algebra of type $A_n$.
\end{abstract}

\maketitle


\section{Introduction}
The Orlik-Terao algebra $OT_n$ is the subalgebra of rational functions on $\C^n$ generated by $\frac{1}{x_i-x_j}$ for all $i\neq j$.
It has been intensively studied in \cite{Terao02,PS06,ST09,Berget, Schenck,DenhamGT,Le14,Liu16,EPW16,MPY16,MMPR21}. 
Only recently, has an attempt to describe the symmetric group action on $OT_n$  been made by Moseley, Proudfoot, and Young \cite{MPY16}.
They provided a recursive algorithm for computing the graded Frobenius character of the $OT_n$.
That algorithm is based on a surprising relation between the Orlik-Terao algebra and the intersection cohomology ring $M_n$ of a certain hypertoric variety constructed from the root system of type $A_n$ \cite{BradenProudfoot, McBreenProudfoot}.

Computation of $M_n$ using the aforementioned algorithm has suggested the following conjecture. 
Let $D_n$ be the cohomology algebra of the configuration spaces of $n$ ordered points in $SU_2$ up to translations.
\begin{conjecture}[{\cite[Conjecture 2.10]{MPY16}}]
For each $n$, there exists an isomorphism of graded $S_n$-representations $D_n \simeq M_n$.
\end{conjecture}
It has been verified for $n\leq 10$ in \cite{MPY16} and for $n \leq 22$ in \cite{MMPR21}.

The algebra $D_n$ has an independent interest, indeed each graded piece is the Whitehouse lift of Eulerian $S_n$-representation up to a sign ($D_n^{k} = \sgn_n \otimes F_{n}^{(n-1-k)}$ see \cite{GS87,Hanlon,Whitehouse, EarlyReiner}).
The Eulerian representations appear also in the study of the free Lie algebra \cite{Reutenauer}.
These representations are used to decompose the Hochshild Cohomology and Cyclic Cohomology in simpler pieces \cite{Whitehouse}.
Moreover, $D_n$ appears in the Hochschild-Pirashvili homology of a wedge of circles and in the weight-zero compactly supported cohomology of $\mathcal{M}_{2,n}$ \cite{GadishHainaut}.

Some tentatives to prove the Moseley-Proudfoot-Young conjecture failed for two reason:
firstly the only known formula describing $D_n$ is
\[ C_n= (V_{(n)} \oplus q V_{(n-1,1)}) \otimes D_n,\]
where $V_\lambda$ is the Schur representation and $C_n$ is the cohomology of the configuration space of $\R^3$.
Although there is an explicit formula for $C_n$ involving plethysm (\Cref{ch_C_n}), inverting the Kronecker (tensor) product is very difficult.
The second issue is that the recursive formula of \cite{MPY16} for $M_n$ is complicate and involves plethysm, Kronecker product and the character of $C_n$. 

We overcome the first problem providing a new virtual formula for the graded Frobenius character of $D_n$ (\Cref{thm:ch_D_n}) by using the Cohen–Taylor-Totaro-K\v{r}iz spectral sequence \cite{CT77,Totaro,Kriz}.
Instead of working on the recursive formula \cite[Theorem 3.2]{MPY16}, we use the isomorphism of graded $S_n$-representations
\[ OT_n \simeq M_n \otimes R_n\]
provided in \cite[Proposition 7]{PS06}, where $R_n$ is the symmetric algebra on $V_{n-1,1}$.
Then we virtually invert $R_n$ (\Cref{lemma:cancelation}) with respect the Kronecker product and we prove the conjecture by induction on $n$ (\Cref{thm:main}) relying on a certain subspace $T_n$ of $OT_n$ (\Cref{lemma:T_Lie}).
Finally, we obtain an explicit formula for the character of $OT_n$ (\Cref{cor:ch_OT_n}) and the generating functions for the characters of $D_n$ and of $OT_n$ (\Cref{cor:gen_funct}).

\section{definitions}
We introduce the main objects of study and some notations.
The Orlik-Terao algebra was introduced in \cite{Terao02}, in type $A_{n-1}$ the definition specializes as follow.

\begin{definition}
\label{def:OT}
The Orlik-Terao algebra of type $A_{n-1}$ is the ring $OT_n= \Q[e_{ij}]/I^{OT}_n$ generated by $e_{i,j}$ for distinct $i,j \in [n]$ and relations  $I_n^{OT}$ given by:
\begin{itemize}
\item $e_{ij} + e_{ji}=0$ for all $i, j$ distinct,
\item $e_{ij} e_{jk} + e_{jk}e_{ki} + e_{ki}e_{ij}=0$ for all $i, j, k$ distinct.
\end{itemize}
\end{definition}

\begin{definition}
Let $C_n^\bullet := H^{2 \bullet}(\Conf_n(\R^3); \Q)$ be the cohomology algebra of the ordered configuration space of $n$ points in $\R^3$.
\end{definition}

The ring $C_n$ can be presented as quotient of $OT_n$ by the equations
\begin{itemize}
\item $e_{ij}^2=0$ for all $i, j$ distinct.
\end{itemize}
The above presentation was studied for the first time in \cite{OrlikTerao}.

\begin{definition}
Let $D_n^\bullet := H^{2\bullet}(\Conf_n(SU_2)/SU_2; \Q)$ be the cohomology algebra of the ordered configuration space of $n$ points in $SU_2$ up to translations.
\end{definition}

The algebra $D_n$ can be presented as $\Q[e_{ij}]/I^{D}_n$ generated by $e_{i,j}$ for distinct $i,j \in [n]$ and relations  $I_n^{D}$ given by:
\begin{itemize}
\item $e_{ij} + e_{ji}=0$ for all $i, j$ distinct,
\item $(e_{ij} + e_{jk} + e_{ki})^2=0$ for all $i, j, k$ distinct,
\item $\sum_{j\neq i} e_{ij}=0$ for all $i \in [n]$.
\end{itemize}
This presentation is due to Matherne, Miyata, Proudfoot, and Ramos \cite[Theorem A4]{MMPR21}.

\begin{definition}
Let $M_n= OT_n/I_n^M$ be the quotient of the Orlik-Terao algebra by the relations:
\begin{itemize}
\item $\sum_{j\neq i} e_{ij}=0$ for all $i \in [n]$.
\end{itemize}
\end{definition}

The algebra $M_n$ was originally defined in a geometric way.
\begin{theorem}[{\cite[Theorem A.6.]{MMPR21}}]
The algebra $M_n^\bullet$ is isomorphic to $\operatorname{IH}^{2 \bullet}(X_n;\Q)$, the intersection cohomology of a hypertoric variety $X_n$ associated with the root system of the Lie
algebra $\mathfrak{sl}_n$.
\end{theorem}

We use the standard notation for symmetric polynomial: let $h_\lambda$, $e_\lambda$, $s_\lambda$, $p_\lambda$ for $\lambda \vdash n$ a partition of $n$ be the complete homogeneous, elementary, Schur, and power sum symmetric polynomials, respectively.
Given a graded $S_n$-representation $V$ we consider the graded Frobenius character $\ch_V(q)$, frequently will omit the dependence on $q$.
As an example if $V_\lambda$ is the irreducible Schur representation in degree zero, then $\ch_{V_{\lambda}}=s_\lambda$.

We denote the \emph{plethysm} of symmetric functions $f,g$ by $f[g]$.  
For $W$ a representation of $S_j$ we denote $\widetilde{W}= W^{\boxtimes m}$ the representation of the wreath product $S_j \wr S_m = (S_j)^{\times m} \rtimes S_m$, where $S_j^{\times m}$ acts coordinatewise and $S_m$ by permuting the coordinates.
Let $V$ be a representation of $S_m$ and $V \otimes \widetilde{W}$ be the representation of $S_j \wr S_m$ where $S_j^{\times m}$ acts only on $\widetilde{W}$ and $S_m$ on both factors.
The group $S_j \wr S_m$ is naturally a subgroup of $S_{jm}$, the main property of the plethysm is 
\[\ch_{\Ind_{S_j \wr S_m}^{S_{jm}} V \otimes \widetilde{W}}= \ch_V [\ch_W]. \]

Let $Lie_n$ be the submodule of the multilinear part of the free Lie algebra on $n$ generators. 
As $S_n$ representation $Lie_n= \Ind_{Z_n}^{S_n} \zeta_n$ where $Z_n$ is the cyclic group generated by an $n$-cycle in $S_n$ and $\zeta_n$ is a primitive root of the unity.
We denote by $l_j$ its character, cf.\ \Cref{remark:final} for an explicit description.
The following result is due to Sundaram and Welker \cite[Theorem 4.4(iii)]{SundaramWelker}, see also \cite[Theorem 2.7]{HershReiner}.
\begin{proposition}
\label{ch_C_n}
The graded character of $C_n$ is
\[ \ch_{C_n} = \sum_{\lambda \vdash n} q^{n-\ell(\lambda)} \prod_{j \geq 1} h_{m_j}[l_j],\]
where $\lambda = (1^{m_1},2^{m_2}, \dots, n^{m_n})$ in the exponential notation and $\ell(\lambda)=\sum_j m_j$ is the number of blocks.
\end{proposition}

Finally, we define $R_n= S^\bullet V_{(n-1,1)}$ and $\Lambda_n = \Lambda^\bullet V_{(n-1,1)}$ be the symmetric (resp.\ alternating) algebra on the standard representation of $S_n$.
We regard $V_{(n-1,1)}$ in degree one, hence $\ch_{\Lambda_n}= \sum_{i=0}^{n-1} q^i s_{n-i,1^i}$.
See \Cref{remark:final} for an expression of $\ch_{R_n}$ in term of Schur polynomials.

\section{Graded Frobenius character of \texorpdfstring{$D_n$}{D}}

In this section we provide a virtual formula for $\ch_{D_n}$ that will be used in the proof of \Cref{thm:main}.
We denote by $\ch'_V$ the expression $\ch_V(-q)$ for $V$ a graded $S_n$-representation.
Let $P_n$ be the $S_n$-representation by permutations, i.e.\ $P_n= V_{(n-1,1)} \oplus V_{(n)}$.
For a partition $\lambda = (1^{m_1},2^{m_2}, \dots, n^{m_n})  \vdash n$ let $S_\lambda$ be the subgroup of $S_n$ stabilizing $\lambda$, i.e. $S_\lambda= \prod_{j\geq 1} S_j \wr S_{m_j}$.

\begin{theorem}
\label{thm:ch_D_n}
The graded character of $D_n$ is:
\begin{equation}
\label{eq:ch_D_n}
\ch_{D_n}(q)= \sum_{\lambda \vdash n}   \frac{q^{n-\ell(\lambda)}}{1-q} \prod_{j \geq 1}  \ch'_{\Lambda^\bullet (P_{m_j})}[l_j].
\end{equation}
\end{theorem}
\begin{proof}
We consider the Cohen–Taylor-Totaro-K\v{r}iz spectral sequence $E_\bullet(SU_2)$ \cite{CT77,Totaro,Kriz} that converge to $H^{\bullet}(\Conf_n(SU_2))$.
In our case since $SU_2$ is $3$-dimensional and has nonzero cohomology only in degree $0$ and $3$, we have that $E_2^{p,q}=0$ if $3 \nmid p$ and $2 \nmid q$.
The $S_n$-representation on the second page is described in \cite[Theorem 3.15]{AAB}:
\begin{equation}
\label{eq:E_2}
 E_2^{3p,2q} = \bigoplus_{\substack{\lambda \vdash n \\ \ell(\lambda)=n-q}} \Ind_{S_\lambda}^{S_n}  \left( \boxtimes_{j} (\Ind_{Z_{j}}^{S_{j}} \zeta_{j})^{\boxtimes m_j} \otimes \Res_{W_\lambda}^{S_{\ell(\lambda)}} \Lambda^p P_{\ell(\lambda)} \right).
\end{equation}
Since $\Res_{W_\lambda}^{S_{\ell(\lambda)}} P_{\ell(\lambda)} = \oplus_{j\geq 1} P_{m_j}$ we have
\begin{equation}
\label{eq:ch_E_2}
\ch_{E_2}(s,t)= \sum_{\lambda \vdash n} t^{2(n-\ell(\lambda))} \prod_{j \geq 1} \ch_{\Lambda^\bullet P_{m_j}}(s^3)[l_j].
\end{equation}
Topologically $SU_2 \simeq S^3$ is a formal orientable manifold, the only nonzero differential of $E_\bullet (SU_2)$ is $\dd_3$ as observed in \cite[\S 1.10]{Petersen} and in \cite[Section 2]{Getzler99}.
The differential $\dd_3$ is compatible with the $S_n$-action by the functoriality property of the spectral sequence.
It follows 
\begin{equation}
\label{eq:E_2_E_inf}
\ch_{E_2}(-q^2,q^3)=\ch_{E_\infty}(-q^2,q^3),
\end{equation}
because this is the right evaluation that simplifies the coimage of $\dd_3$ with its image.

Consider the map $f \colon (\mathbb{R}^3)^{n-1} \to (SU_2)^n$ defined by $(x_1, \dots, x_{n-1}) \mapsto (x_1, \dots, x_{n-1},e)$ where $e$ is the identity of $SU_2$ and $\mathbb{R}^3$ is identified with $SU_2 \setminus \{e\}$.
It has a retraction defined by 
\[(g_1, g_2, \dots, g_n) \mapsto (g_n^{-1}g_1,  g_n^{-1}g_2, \dots, g_n^{-1}g_{n-1}).\]
Both maps restrict to the subspaces $\Conf_{n-1}(\mathbb{R}^3) \to \Conf_n (SU_2)$ and this implies that $E_\bullet (\mathbb{R}^3)$ is a direct addendum of $E_\bullet(SU_2)$.
Notice that $\Conf_{n-1}(\mathbb{R}^3) \times SU_2 \simeq \Conf_n(SU_2)$ via the map $((x_1, \dots, x_{n-1}),g) \mapsto g \cdot f(\underline{x})$, hence $E_\infty(SU_2)= E_\infty(\mathbb{R}^3) \otimes H^\bullet(SU_2)$ as graded vector spaces.
Since $E_2(\mathbb{R}^n)$ is supported on the column $p=0$, so is $E_\infty(\mathbb{R}^n)$.
Therefore $E_\infty (SU_2)$ is supported only on the column $p=0$ and $p=3$, indeed the even cohomology of $\Conf_n(SU_2)$ is supported in degrees $(0, 2q)$ and the odd one in degrees $(3,2q)$.
So 
\[\ch_{E_\infty}(s,t)= \ch_{H^{\textnormal{even}}(\Conf_n(SU_2))}(t)+ s^3 t^{-3} \ch_{H^{\textnormal{odd}}(\Conf_n(SU_2))}(t).\]
Let $\pi \colon \Conf_n(SU_2) \to \Conf_n(SU_2)/SU_2$ be the natural projection, it is a $S_n$-equivariant fiber bundle. 
The Leray-Hirsch theorem for rational cohomology asserts that $H(\Conf_n(SU_2))$ is a free $H( \Conf_n(SU_2)/SU_2)$-module with basis given by $1, \omega$ for any nonzero $\omega \in H^3(\Conf_n(SU_2))$.
The module structure is given by $\pi^*$ so it is $S_n$-equivariant. 
We observe that $S_n$ acts trivially on $H^0(\Conf_n(SU_2))$ and on $H^3(\Conf_n(SU_2))$, because the latter is a $1$-dimensional quotient of $E^{3,0}_2(SU_2) \cong P_n$.
Therefore
\begin{align*}
\ch_{H^{\textnormal{even}}(\Conf_n(SU_2))}(t) &= \ch_{H(\Conf_n(SU_2)/SU_2)}(t) = \ch_{D_n}(t^2), \\
\ch_{H^{\textnormal{odd}}(\Conf_n(SU_2))}(t) &= t^3 \ch_{H(\Conf_n(SU_2)/SU_2)}(t) = t^3 \ch_{D_n}(t^2).
\end{align*}
We have $\ch_{E_\infty}(s,t)= (1+s^3) \ch_{D_n}(t^2)$ and together with eq.\ \eqref{eq:ch_E_2} and \eqref{eq:E_2_E_inf} they imply
\[ (1-q^6) \ch_{D_n}(q^6) =  \sum_{\lambda \vdash n} q^{6(n-\ell(\lambda))} \prod_{j \geq 1} \ch_{\Lambda^\bullet P_{m_j}}(-q^6)[l_j]. \]
That is our claim.
\end{proof}

\begin{remark}
The formula \eqref{eq:ch_D_n} has $(1-q)$ in the denominator and seems to be an infinite series.
However it can be written as a polynomial in $q$ of degree $n-1$:
\[ \ch_{D_n}(q)= \sum_{\lambda \vdash n}   q^{n-\ell(\lambda)} (1-q)^{c_\lambda -1} \prod_{j \geq 1}  \ch'_{\Lambda_{m_j}}[l_j],\]
where $c_\lambda=\lvert \{ j \mid m_j \neq 0 \}\rvert$.
Furthermore, since the left hand side is a polynomial in $q$ of degree $n-2$, the coefficient  of $q^{n-1}$ in the right hand side must be zero.
\end{remark}

\section{Proof of the MPY conjecture}
Now we prove the conjecture and provide a new formula for the character of the Orlik-Terao algebra.
The Kronecker product of two symmetric function $f*g$ is the linear extension of the tensor product for representation, i.e. $\ch_{V\otimes W}= \ch_V *\ch_W$.

\begin{theorem}[{\cite[Proposition 7]{PS06}}]
\label{thm:OT_is_M_x_R}
For each $n$ the equation
\begin{equation*}
\ch_{OT_n}= \ch_{M_n} * \ch_{R_n} 
\end{equation*}
holds.
\end{theorem}

\begin{lemma}
\label{lemma:cancelation}
Let $V$ be any representation of the symmetric group $S_n$. We have:
\[ \ch_{S^\bullet V} * \ch'_{\Lambda^\bullet V} = s_n. \]
\end{lemma}
\begin{proof}
The Koszul complex for the ring $S^\bullet V$ is a free resolution of $\Q= S^\bullet V/(V)$. The bigraded character of the Koszul complex is $\ch_{S^\bullet V}(s)*\ch_{\Lambda^\bullet V}(t)$, hence by exactness we have $\ch_{S^\bullet V}(q)*\ch_{\Lambda^\bullet V}(-q)=s_{n}$.
\end{proof}
It follows that $\ch_{R_n}$ is invertible with respect to the Kronecker product, whose inverse is $\ch'_{\Lambda_n}$. 

\begin{lemma}
\label{lemma:plethysm}
Let $g$ be a symmetric function of degree $j$ and $m$ a positive integer. We have
\[
\ch'_{\Lambda^\bullet P_{m}} [g] = h_{m}[(1-q) g].
\]
\end{lemma}
\begin{proof}
Using the identity $h_{n-k} e_k = s_{n-k,1^k}+ s_{n-k+1,1^{k-1}}$ we obtain
\[ \ch'_{\Lambda^\bullet P_{n}}= (1-q) \sum_{k=0}^{n-1} (-q)^k s_{n-k,1^k}=  \sum_{k=0}^{n} (-q)^k h_{n-k} e_k. \] 
Recall the subtraction formula (see for example in \cite[\S 3.3]{LoehrRemmel})
\[ h_m[f-g] = \sum_{i=0}^m (-1)^{k} h_{m-k} [f]e_{k} [g], \]
we obtain
\begin{align*}
h_{m}[(1-q) g] &=\sum_{k=0}^m (-1)^{k} h_{m-k} [g]e_{k} [q g]\\
&=\sum_{k=0}^m (-q)^{k} (h_{m-k}e_{k})[g]\\
&=  \ch'_{\Lambda^\bullet P_{m}} [g]. \qedhere
\end{align*}
\end{proof}

Using the Lemma above we can rewrite the character of $D_n$ as follow.
\begin{corollary}
\label{cor:ch_D_2}
The graded character of $D_n$ is
\begin{equation}
\label{eq:ch_D_2}
\ch_{D_n}(q)=  \frac{1}{1-q} \sum_{\lambda \vdash n}   \prod_{j \geq 1}  h_{m_j}[q^{j-1} (1-q)l_j].
\end{equation}
\end{corollary}
\begin{proof}
It follows from \Cref{thm:ch_D_n} and \Cref{lemma:plethysm}.
\end{proof}

\begin{lemma}
\label{lemma:tensor_induction}
Let $\lambda =(1^{m_1}, 2^{m_2}, \dots)$ be a partition of $n$ and $g_j, f_{m_j}$ be symmetric functions of degree $j$ and $m_j$ respectively.
We have:
\[ \ch'_{\Lambda^\bullet P_n}* \prod_{j\geq 1}  f_{m_j}[g_j] =  \prod_{j\geq 1}  f_{m_j} [g_j * \ch'_{\Lambda^\bullet P_j} ].\]
\end{lemma}
\begin{proof}
Firstly observe that 
\[ \Res^{S_n}_{\prod_{j\geq 1} S_{jm_j}} P_n = \bigoplus_{j \geq 1} P_{jm_j},\]
and so
\[\Res^{S_n}_{\prod_{j\geq 1} S_{jm_j}} \Lambda^\bullet P_n= \bigotimes_{j\geq 1} \Lambda^\bullet P_{jm_j}. \]
Using the projection formula (sometimes called Frobenius reciprocity) we obtain:
\[ \ch'_{\Lambda^\bullet P_n}* \prod_{j\geq 1}  f_{m_j}[g_j] = \prod_{j\geq 1}   \ch'_{\Lambda^\bullet P_{jm_j}}* f_{m_j}[g_j].\]
Thus it is enough to show 
\[\ch'_{\Lambda^\bullet P_{jm}}* f[g]= f [g * \ch'_{\Lambda^\bullet P_j} ].\] 
This last equality is linear and multiplicative in the entry $f$: the linearity is trivial and the multiplicativity follow from the argument above
\begin{align*}
 \ch'_{\Lambda^\bullet P_{jm}}* (f_1f_2)[g] &= \ch'_{\Lambda^\bullet P_{jm}}* (f_1[g] f_2[g]) \\
 &= (\ch'_{\Lambda^\bullet P_{jm_1}}* f_1[g] )(\ch'_{\Lambda^\bullet P_{jm_2}}* f_2[g] ).
\end{align*}
Therefore we may assume $f=p_m$. 
Again $\ch'_{\Lambda^\bullet P_{jm}}* p_m[g]= p_m [g * \ch'_{\Lambda^\bullet P_j} ]$ is linear and multiplicative in the entry $g$ and so we reduce to the case $g=p_j$.

It remains to prove that $\ch'_{\Lambda^\bullet P_{jm}}* p_{jm}= p_m [p_j * \ch'_{\Lambda^\bullet P_j} ]$.
Since $(p_\lambda)_{\lambda}$ are orthogonal idempotent with respect to the Kronecker product 
\[\ch'_{\Lambda^\bullet P_{n}}* p_{n} = \chi'_{\Lambda^\bullet P_{n}} (c_n) p_n \]
where $\chi'_V(\sigma)$ is the graded character of $\sigma \in S_n$ with $q$ replaced by $-q$ and $c_n \in S_n$ be an $n$-cycle.
It is easy to see that 
\[ \chi'_{\Lambda^\bullet P_{n}} (c_n)=1+(-1)^{n-1}(-q)^n=1-q^n\]
 on the canonical base of $\Lambda^\bullet P_{n}$: let $(v_i)_i$ the standard base of $P_n$,  the product of some $v_j$ is invariant for $c_n$ if and only if each generator appears a fixed number of times (i.e.\ $0$ or $1$ times).
Finally the equalities
\begin{align*}
p_m [p_j * \ch'_{\Lambda^\bullet P_j} ]
&= p_m[(1-q^j)p_j] \\
&= (1-q^{jm}) p_{jm} \\
&= \ch'_{\Lambda^\bullet P_{jm}}* p_{jm}
\end{align*}
conclude the proof.
\end{proof}

For each monomial $m= \prod_{k} e_{i_k,j_k} \in \Q[e_{i,j}]$ we define the \emph{support} of $m$ as the finest set partition $B(m) \vdash [n]$ such that for all $k$ $i_k$ and $j_k$ belong to the same block of $B(m)$.
We also define the \emph{type} of $m$ as the partition $\lambda(m) \vdash n$ collecting the size of blocks of $B(m)$.
Notice that the relations defining $OT_n$ (\Cref{def:OT}) are sum of monomials with the same support, hence the notion of support and type are well defined in $OT_n$. 
Moreover, monomials with different supports are linearly independent.

For $B\vdash [n]$ a set partition let $T_B \subset OT_n$ be the vector space generated by all monomials $m$ such that $B(m)=B$.
For $S\subseteq [n]$ we define $T_S=T_{B}$ where $B$ is the finest set partition of $[n]$ with a block equal to $S$.
Given two monomials $m,n$ such that $mn \neq 0$ in $OT_n$, we have that $B(mn)$ is the finest set partition coarsening both $B(m)$ and $B(n)$, hence 
\[T_B \cong \bigotimes_{i=1}^{\lvert B \rvert} T_{B_i}.\]

Consider a partition $\lambda \vdash n$, let $T_\lambda$ be the vector space generated by all monomials of type $\lambda$.
Choose a set partition $B_\lambda \vdash [n]$ whose blocks $B_i$ are of length $\lambda_i$ and let $S_{B_\lambda}$ be the subgroup of $S_n$ stabilizing $B_\lambda$, if $\lambda=(1^{m_1}, 2^{m_2}, \dots, n^{m_n})$ then $S_{B_\lambda} \cong \prod_{j\geq 1} S_j \wr S_{m_j}$.
We have 
\[T_\lambda \cong \Ind_{S_{B_\lambda}}^{S_n} T_{B_\lambda}\]
 as representation of $S_n$, where $S_{B_i}$ acts on the factor $T_{B_i}$ of $T_{B_\lambda}= \otimes_{i=1}^{\lvert B \rvert} T_{B_i}$ and $S_{m_j}$ permutes the $m_j$ factors of size $j$.
For the sake of notation we set $T_n = T_{(n)}$.

\begin{lemma}
\label{lemma:T_Lie}
We have
\[ \ch_{OT_n}= \sum_{\lambda \vdash n} \prod_{j \geq 1} h_{m_j}[\ch_{T_n}]. \]
\end{lemma}
\begin{proof}
The Orlik-Terao algebra decomposes
\begin{align*}
OT_n =& \bigoplus_{B \vdash [n]} T_{B} \\
=& \bigoplus_{B \vdash [n]} \bigotimes_{i=1}^{\lvert B \rvert} T_{B_i} \\
=& \bigoplus_{\lambda \vdash n} \Ind_{S_{B_\lambda}}^{S_n} \bigotimes_{i=1}^{ \ell(\lambda)} T_{B_i}  \\
=& \bigoplus_{\lambda \vdash n} \Ind_{\prod_j S_{jm_j}}^{S_n} \left( \bigotimes_{j \geq 1} \Ind_{S_j \wr S_{m_j}}^{S_{jm_j}} \widetilde{T_{j}} \right)
\end{align*}
as $S_n$-representation.
Taking the character we obtain the claimed relation.
\end{proof}

\begin{theorem}
\label{thm:main}
We have
\[\ch_{D_n} = \ch_{M_n}  \]
and
\[ \ch_{T_n}= q^{n-1} l_n * \ch_{R_n}.\]
\end{theorem}
\begin{proof}
We prove both equality by induction on $n$.
The base case $n=1$ is trivial.
For the inductive step we consider:
\begin{equation*}
\begin{aligned}
\ch_{M_n}=& \ch_{OT_n} * \ch'_{\Lambda_n} \\
=&\frac{1}{(1-q)}  \sum_{\lambda \vdash n} \prod_{j \geq 1} h_{m_j}[\ch_{T_j}*\ch'_{\Lambda^\bullet P_j}] \\
=& \ch_{T_n}* \ch'_{\Lambda_n} +\frac{1}{(1-q)}  \sum_{\substack{\lambda \vdash n \\ \lambda \neq (n)}} \prod_{j \geq 1} h_{m_j}[ q^{j-1}(1-q) l_j].
\end{aligned}
\end{equation*}
The first equality follows from \Cref{thm:OT_is_M_x_R} and \Cref{lemma:cancelation}.
The second one follows from \Cref{lemma:T_Lie} and \Cref{lemma:tensor_induction} together with the identity $\ch'_{\Lambda^\bullet P_j}= (1-q) \ch'_{\Lambda_j}$.
The last one follows from the inductive hypothesis and \Cref{lemma:cancelation}.
We have proven the identity
\[ \ch_{M_n}- \ch_{T_n}*\ch'_{\Lambda_n} = \frac{1}{(1-q)}  \sum_{\substack{\lambda \vdash n \\ \lambda \neq (n)}}   \prod_{j \geq 1} h_{m_j}[ q^{j-1} (1-q) l_j] = \ch_{D_n}- q^{n-1} l_n,\]
where the last equality is given by \Cref{cor:ch_D_2}.
Since $\ch_{D_n}$ and $\ch_{M_n}$ has degree less than $n-1$ and $ \ch_{T_n}*\ch'_{\Lambda_n}$ bigger than $n-2$,  $\ch_{M_n}=\ch_{D_n}$ and $\ch_{T_n}*\ch'_{\Lambda_n} = q^{n-1} l_n$ hold.
Therefore $\ch_{T_n}= q^{n-1} l_n * \ch_{R_n}$.
\end{proof}

\begin{corollary}
\label{cor:ch_OT_n}
We obtain the character of $OT_n$:
\begin{equation}
\label{eq:OT_n}
\ch_{OT_n}= \sum_{\lambda \vdash n} q^{n-\ell(\lambda)}\prod_{j\geq 1} h_{m_j}[l_j * \ch_{R_j}].
\end{equation} 
\end{corollary}
\begin{proof}
It follows from \Cref{thm:main} and \Cref{lemma:T_Lie}.
\end{proof}

An important object for the proof of \Cref{thm:main} is the $R_n$-module $T_n$. It is a submodule of the free module $OT_n$ and its Frobenius character is equal to the one of the free module $R_n \otimes_\Q T_n^{n-1}$.
This observations lead to the following conjecture:
\begin{conjecture}
The $R_n$-module $T_n$ is free.
\end{conjecture}

\begin{remark}
\label{remark:final}
The formula \eqref{eq:OT_n} is completely explicit because $\ch_{R_j}$ and $l_j$ are known.
Indeed
\[ \ch_{R_n}=  (1-q) \sum_{\lambda \vdash n} s_{\lambda}(1,q,q^2,... ) s_\lambda = (1-q) h_n \left[\frac{X}{1-q} \right]\]
by \cite[Section 5.6]{Procesi03} or \cite[Exercise 7.73]{StanleyVol2} where $X=h_1$.
Moreover,
\[ l_n = \frac{1}{n} \sum_{d \mid n} \mu(d) p_d^{\frac{n}{d}},\]
by \cite[Theorem 8.3]{Reutenauer}, $l_n$ is known as the Lyndon symmetric function or as Gessel-Reutenauer symmetric function \cite{GesselReutenauer}.
\end{remark}

Let $\Exp$ be the plethystic exponential defined by 
\[\Exp(f)= \operatorname{exp} \left( \sum_{k\geq 1} \frac{p_k[f]}{k} \right) = \sum_{k\geq 0} h_k[f],
\]
and $\operatorname{Log}$ be its inverse.
We define the symmetric functions 
\[L=\sum_{n\geq 1} q^{n-1}t^n l_n = - \frac{ \operatorname{Log}(1-qtX)}{q}\]
and $H=\sum_{k\geq 0} h_k$.

\begin{corollary}
\label{cor:gen_funct}
The generating functions for $\ch_{D}$ and $\ch_{OT}$ are:
\begin{gather}
\sum_{t \geq 1} \ch_{D_n}(q)t^n = \frac{1}{1-q} (\Exp((1-q)L)-1), \\
\sum_{t \geq 1} \ch_{OT_n}(q)t^n =  \Exp\left( (1-q)L*H \left[\frac{X}{1-q} \right] \right)-1.
\end{gather}
\end{corollary}
\begin{proof}
Let $f$ be a symmetric function and call $f_j$ be the homogeneous part of degree $j$.
Assume that $f$ has zero constant term, i.e. $f= \sum_{j\geq 1} f_j$, then
\begin{align*}
\Exp(f) &= \prod_{j\geq 1} \Exp(f_j) \\
&=    \prod_{j\geq 1} \sum_{m\geq 0} h_m[f_j] \\
&= \sum_{\lambda} \prod_{j\geq 1} h_{m_j}[f_j],
\end{align*}
where the sum is taken over all partitions $\lambda=(1^{m_1}, 2^{m_2}, \dots)$.
The corollary follows by taking $f=(1-q)L$ and $f=(1-q)L*H[(1-q)^{-1}X]$.
\end{proof}

Formulas of this paper are checked and implemented in SageMath \cite{Sage}.
The code is available at 
\begin{center}
\url{https://github.com/paga92/character_OT}.
\end{center}

\section*{acknowledgement}
I would like to thank Nir Gadish for useful discussion and for pointing out an error in the first version  and Alessando Iraci for suggested a simplified proof of \Cref{lemma:plethysm}.
Part of the work was carried out during my stay at the institute for Computational and Experimental Research in Mathematics in Providence, RI, during the Braids in Representation Theory and Algebraic Combinatorics program.

\bibliographystyle{amsalpha}
\bibliography{bib}

\end{document}